\documentclass[12pt]{article}
\usepackage[cp1251]{inputenc}
\usepackage{amsfonts,enumerate,amsmath,amssymb}

\newtheorem{theorem}{Theorem}

\def\C{{\mathbb C}}
\def\R{{\mathbb R}}
\def\Z{{\mathbb Z}}
\def\Q{{\mathbb Q}}

\begin{document}

\title{A canonical basis of two-cycles on a $K3$ surface}
\author{I.A. Taimanov
\thanks{Sobolev Institute of Mathematics, Academician Koptyug avenue 4, 630090, Novosibirsk, Russia, and Department of Mathematics and Mechanics, Novosibirsk State University, Pirogov street 2, 630090 Novosibirsk, Russia; 
e-mail: taimanov@math.nsc.ru. \newline The work was supported by RSF (grant 14-11-00441).}}
\date{}
\maketitle

In this article we obtain an explicit form for a canonical basis of two-cycles on a $K3$ surface. This basis is realized by formal sums of smooth submanifolds.

By a canonical basis we mean a basis in which the intersection form
$$
H_2(X;\Z) \cap H_2(X;\Z) \to \Z = H_0(X;\Z)
$$
looks as follows
\begin{equation}
\label{form}
E_8(-1) \oplus E_8(-1) \oplus H \oplus H \oplus H,
\end{equation}
where
\begin{equation}
\label{e8}
E_8(-1) = (-1) \cdot
\left(
\begin{array}{cccccccc}
2 & 0 & -1 & 0 & 0 & 0 & 0 & 0 \\
0 & 2 & 0 & -1 & 0 & 0 & 0 & 0 \\
-1 & 0 & 2 & -1 & 0 & 0 & 0 & 0 \\
0 & -1 & -1 & 2 & -1 & 0 & 0 & 0 \\
0 & 0 & 0 & -1 & 2 & -1 & 0 & 0 \\
0 & 0 & 0 & 0 & -1 & 2 & -1 & 0 \\
0 & 0 & 0 & 0 & 0 & -1 & 2 & -1 \\
0 & 0 & 0 & 0 & 0 & 0 & -1 & 2
\end{array}
\right)
\end{equation}
and
\begin{equation}
\label{h}
H = \left(\begin{array}{cc} 0 & 1 \\ 1 & 0 \end{array}\right).
\end{equation}

Originally the intersection form of a $K3$ surface was found in \cite{Milnor} by the general reasonings:
according to the Milnor theorem proved in the same article (see also \cite{Serre}),
an even indefinite unimodular form over $\Z$ is uniquely defined by its rank $r$ and
signature $\tau$ and is equal to
$$
\left(\frac{-\tau}{8}\right) E_8(-1) \oplus \left(\frac{r+\tau}{2}\right) H.
$$
Since in our case the rank equals to the second Betti number
$b_2(X)=22$, which is derived by using only rational cohomology, the signature $\tau=-16$ is derived from
Hodge theory, and the form is even because the second Stiefel--Whitney class
$w_2(X)$ vanishes, the intersection form of a $K3$ surface $X$ is given by (\ref{form}).
Detailed proofs in the same spirit are given, for instance, in \cite{AS,BHPV}.
Until recently a canonical basis in $H_2(X;\Z)$ was not found explicitly and, in particular,
in the recent lectures on $K3$ surfaces there was expressed a wish to fill this gap \cite{H}.

Since all $K3$ surfaces are not only pairwise homeomorphic but also dif\-feo\-morphic (see, for instance, \cite{AS}), it is enough for us to
use an exact model, of such a complex surface, for which we take the Kummer manifold.

Let us consider a four-dimensional torus $T^4$, which is the quotient of a two-dimensional complex space
$\C^2$ under the action (by translations) of a lattice $\Lambda$ of rank $4$: $T^4 = \C^2/\Lambda$.
The reflection $\iota: \C^2 \to \C^2, \iota(z)=-z$, generates a $\Z_2$-action on $T^4$. This $\Z_2$--action is not free
and has 16 fixed points which lie in the half-periods of the lattice, i.e. in $\frac{1}{2}\Lambda$. Every fixed point
has a $\iota$-invariant disc neighborhood $D^4$ with a spherical boundary $\partial D^4 = S^3$.
By the projection $T^4 \to T^4/\Z_2$ such a disc boundary is mapped into a cone over the boundary
$\partial D^4/\Z_2 = \R P^3$ and with the apex at the fixed point.

The Kummer manifold $X$  is obtained from the quotient space $T^4/\Z_2$
by the blowup of singularities corresponding to 16 fixed points of $\iota$. For every singular point
the blowup (the $\sigma$-process)  consists in removing of the corresponding cone
$D^4/\Z_2$ and gluing to the emerging boundary hypersurface $\R P^3$ a smooth manifold $M^4$ with the same boundary.
This manifold is fibered by two-discs over the two-sphere and the fibration $M^4 \to S^2$ is fiberwise embedded
into the complex line bundle $V^4 \stackrel{\C}{\longrightarrow}  \C P^1 = S^2$ which is the square $\gamma^2 = O(-2)$ of the standard linear (or tautological) bundle $\gamma = O(-1)$ over $\C P^1$.

By this surgery every singular point $L \in T^4/\Z_2$ is replaced by a two-sphere $S^2_L =\C P^1$ and there is a natural
projection $M^4 \setminus S^2_L \to D^4/\Z_2\setminus L$ which is a diffeomorphism.
The self-intersection index of $S^2_L$ is equal to $(-2)$ for every singular point $L$.
In particular, this implies that this sphere generates a nontrivial homology class which we denote by $[L]$.
The self-intersection index of this cycle can be calculated by using the standard methods of algebraic geometry.
From the topological point of view this is as follows. Let us mean by the zero section of a vector bundle
the section that assigns the zero vector to every point.
It is known that as a complex linear bundle the tangent bundle
$TS^2$ of a two-sphere is isomorphic to $\bar{\gamma}^2$, where the bundle $\bar{\gamma}$ is dual to
$\gamma$. Let us perturb the zero section $\eta_0$ of the tangent bundle $TS^2 \to S^2$ to achieve a section
$\eta$ that intersects $\eta_0$ transversally.  By definition, a section of the tangent bundle is a vector field
$\eta$ on $S^2$. In $TS^2$ the index of intersection of $\eta$ with $\eta_0$ is equal to the self-intersection index of $S^2$, and it is also equal to the index of the vector field $\eta$, i.e. to the Euler characteristic of $S^2$: $\mathrm{ind}\, \eta = 2$. The Euler characteristic of the sphere $S^2=\C P^1$ is equal to $\langle c_1(\bar{\gamma}^2, [\C P^1] \rangle = - \langle c_1(\gamma^2), [\C P^1]\rangle$, and $\langle c_1(\gamma^2), [\C P^1]\rangle$ is the self-intersection index of the $[L]$-cycle:
$$
[L] \cap [L] = -2.
$$

Multidimensional Kummer manifolds are obtained analogously
by desin\-gu\-la\-ri\-zations of the quotient-spaces
$T^{2N}/\Z_2$ which have $2^{2N}$ singular points.
In this case every singular point is replaced by $\C P^{N-1}$.
The Kummer manifolds are simply-connected and their rational homology are generated by classes of the two types:
1) the classes that are pulled back to the homology classes of tori;
2) the classes that come from the copies of $\C P^{N-1}$ which replace singular points
(see \cite{Spanier}).

For $K3$ surfaces we describe these generators in detail.

For simplicity we consider the case when $\Lambda = \Z^4$. We denote by $x_1$, $x_2$, $x_3$, $x_4$
the Euclidean coordinates in $\R^4 = \C^2$ and by $z_1=x_1+ix_2, z_2=x_3+ix_4$ the coordinates in $\C^2$.
Denote by $e_1,e_2,e_3,e_4$ the corresponding basis in $\R^4$ and, also, in $\Lambda$.

By $T_{ij}$ we denote the oriented two-dimensional torus, in $T^4$,  which is generated by the vectors $e_i$ and $e_j$.
We assume that the orientation is defined by the positive frame $(e_i,e_k)$. We denote by $[T_{ij}]$ the corresponding homology class.

The intersection indices of such classes are as follows
\begin{equation}
\label{form1}
[T_{ij}] \cap [T_{kl}] = 2 \varepsilon_{ijkl},
\end{equation}
where $\varepsilon_{ijkl} = 0$, if among $i,j,k,l$ there are two coinciding indices, $1$ for even transpositions
$\left(\begin{array}{cccc} 1 & 2 & 3 & 4 \\ i & j & k & l \end{array}\right)$, and $(-1)$, if such a transposition is odd.

If $T_{ij}$ does not pass through fixed points of $\iota$, then it projects into a torus in $T^4/\Z_2$, which gives rise to
a nontrivial homology class of $X$. We denote these homology classes by the same symbols and note that the intersection formula for them is also given by (\ref{form1}). 

 We denote by $L_1,\dots,L_{16}$ the fixed points of $\iota$ and by
 $[L_1]$, $\dots$, $[L_{16}]$ the corresponding two-dimensional homology classes.
 As we showed above,
 \begin{equation}
 \label{form2}
 [L_i] \cap [L_j] = -2 \, \delta_{ij}.
 \end{equation}
Moreover
 \begin{equation}
 \label{form3}
 [T_{ij}]\cap [L_k] = 0 \ \ \ \mbox{for all $i,j,k$}.
 \end{equation}
 Therefore the cycles $[L_i], i=1,\dots,16$, and $[T_{ij}], 1 \leq j<j \leq 4$, form a basis in the rational
 homology group $H_2(X;\Q) = \Q^{22}$ and the lattice generated by them is a sublattice of index $2^{22}$ in $H_2(X;\Z)$, because the intersection form is not unimodular on this sublattice and has the determinant equal to
 $-2^{22}$.

 To find a canonical basis in $H_2(X;\Z)$ we have to consider another cycles which are linear combinations of
 $[L]$- and $[T]$-cycles over rational numbers.

 Let us consider a two-dimensional torus $T^2 \subset T^4$ that passes through fixed points. Since it passes through a fixed point (a half-period) and is generated by two vectors from the lattice $\Lambda$, it has to pass through four fixed points. Let us denote by $T_{ijkm}$ the torus that passes through $L_i,L_j,L_k,L_m$. It is invariant under $\iota$.
Let us consider the square (in $T_{ijkm}$) with vertices at $L_i$, $L_j$, $L_k$, and $L_m$ which are ordered so as to be successively traversed while circumventing the boundary of the square.
We assume that such a circumventing a positive orientation of the boundary of the square and therewith the agreed orientations of the square and of the torus $T_{ijkm}$.

 Let us apply a surgery to $T_{ijkm}/\Z_2$ to obtain by desingularization smooth manifolds whose intersection indices with
 $L_n$ are equal to $0$ or $1$. We consider the two cases.

{\sc Case 1.} 
Let $T_{ijkm} \subset \C^2/\Z^4$ be the quotient of the complex line $\C = \{z_2=\mathrm{const}\}$ under the action of $\Z^2$.
The desingularization removes from $T_{ijkm}/\Z_2$ the two-dimensional disc which is the cone with the apex at singular point and replace it by the disc that lies in a fiber of the fibration $\gamma^2: V^4 \to \C P^1 = S^2$ and intersects $L_n$, $n \in \{i,j,k,m\}$, at the unique point which is the zero vector of the fiber.
Thus we obtain the smooth manifold which we denote by $S_{ijkm}$. Since both this manifold and $L_n$ are complex curves, their
intersection index is equal to one.

{\sc Case 2.} 
Let the torus $T_{ijkm}$ be not a smooth submanifold and pass through a fixed point $L$.  It intersects the boundary of the disc-neighborhood $D^4(L)$ of $L$ by the curve $\varkappa_0$, whose orientation agrees with the orientation of the torus and which is $\iota$-invariant. Let us construct some homotopy $\varkappa_t, 0 \leq t \leq 1$, of this curve that
 consists in simultaneous rotation around $L$ and dilation with coefficient $\left(1-\frac{t}{2}\right)$ and center at $L$ and such that the oriented contour $\varkappa_1$ lies in the plane $z_2 = \mathrm{const}$, in which it bounds the positively oriented disc $D_0(L)$. The cylinder swept under the homotopy $\varkappa_t$ and the disc $D_0(L)$ form together the disc $D_1(L)$ which is $\iota$-invariant. We remove from $T_{ijkm}$ the disc $T_{ijkm}\cap D^4(L)$ and replace it by $D_1(L)$. It is evident that by small perturbation of the homotopy $\varkappa_t$ this construction can be smoothed to achieve a $\iota$-invariant smooth submanifold. By applying this construction to all intersections of $T_{ijkm}$ with the fixed points of $\iota$ we obtain a $\iota$-invariant smooth submanifold $T^\prime_{ijkm}$. The blow-up of singularities in $T^\prime_{ijkm}/\Z_2$ has the same form as in the first case and we obtain the smooth submanifold $S_{ijkm} \subset X$ whose intersection indices with $L_n$ are equal either one (for $n \in \{i,j,k,m\}$), or zero (otherwise).

Clearly the submanifolds $S_{ijkm}$ are diffeomorphic to the two-sphere $S^2$ and
 $T^\prime_{ijkm}=T^2 \to S_{ijkm}$ are two-sheeted coverings branched at four points.

We denote by $[S_{ijkm}]$ the homology cycles corresponding to $S_{ijkm}$.

Let us proceed directly to the construction of a canonical basis of cycles.

We split the $[L]$-cycles into two groups corresponding to fixed points with $x_4=0$ and to fixed points with $x_4 = \frac{1}{2}$.

Let us consider the first group and enumerate the cycles in it by the following rule
(on the right-hand sides we give the coordinates of the corres\-pon\-ding fixed point):
$$
[L_1] \leftrightarrow \left(0,0,0\right), \ \
[L_2]  \leftrightarrow \left(0,0,\frac{1}{2}\right), \ \
[L_3]  \leftrightarrow \left(\frac{1}{2},0,0\right), \ \
[L_4] \leftrightarrow \left(\frac{1}{2},0,\frac{1}{2}\right),
$$
$$
[L_5] \leftrightarrow \left(\frac{1}{2},\frac{1}{2},0\right), \ \
[L_6]  \leftrightarrow \left(\frac{1}{2},\frac{1}{2},\frac{1}{2}\right), \ \
[L_7]  \leftrightarrow \left(0,\frac{1}{2},0\right), \ \
[L_8] \leftrightarrow \left(0,\frac{1}{2},\frac{1}{2}\right).
$$

By construction, the cycles $[S_{1357}]$, $[S_{2156}]$, $[S_{5643}]$, and $[S_{3487}]$
have non\-ne\-ga\-tive intersection indices with the cycles of the form $L_n$ and, by decomposing them over $\Q$ in terms
of the basis $\{[L_i], [T_{ij}]\}$ and by (\ref{form1}), (\ref{form2}) and (\ref{form3}), we obtain
\begin{equation}
\label{s}
\begin{split}
[S_{1357}] = - \frac{1}{2}\left([L_1]+[L_3]+[L_5]+[L_7]\right)  +\frac{1}{2}[T_{12}],
\\
[S_{2156}] = - \frac{1}{2}\left([L_1]+[L_2]+[L_5]+[L_6]\right)  + \frac{1}{2}([T_{13}] + [T_{23}]),
\\
[S_{5643}] = - \frac{1}{2}\left([L_3]+[L_4]+[L_5]+[L_6]\right)  + \frac{1}{2}[T_{23}],
\\
[S_{3487}] = - \frac{1}{2}\left([L_3]+[L_4]+[L_7]+[L_8]\right)  + \frac{1}{2}([T_{13}]-[T_{23}])
\end{split}
\end{equation}

{\sc 1) The cycles $w_1,\dots,w_8$.}

Let us introduce the following homology cycles
$$
w_1 =  - [S_{1357}]  - [L_3] - [L_5] - [L_7] =  \frac{1}{2}\left([L_1] - [L_3] - [L_5] - [L_7]\right) - \frac{1}{2}[T_{12}],
$$
$$
w_2= - [S_{2156}] =  \frac{1}{2}\left([L_1] + [L_2] + [L_5] + [L_6]\right) - \frac{1}{2}([T_{13}] + [T_{23}]),
$$
$$
w_3 = ([T_{13}] + [T_{23}]) - [S_{2156}] - [L_1] - [L_2] = 
$$
$$
 = \frac{1}{2}\left( - [L_1] - [L_2] + [L_5] + [L_6]\right) 
+ \frac{1}{2}([T_{13}] + [T_{23}]),
$$
$$
w_4 = - [L_6],
$$
$$
w_5 = - [S_{5643}] - [L_5] =  \frac{1}{2}\left( [L_3] + [L_4] - [L_5] + [L_6]\right) -  \frac{1}{2}[T_{23}],  
$$
$$
w_6 = - [L_4], 
$$
$$
w_7 =  - [S_{3487}] - [L_5] =  \frac{1}{2}\left( - [L_3] + [L_4] + [L_7] + [L_8]\right) - \frac{1}{2}([T_{13}]-[T_{23}]), 
$$
$$
w_8 = - [L_8].
$$

{\sc 2) The cycles $w_9,\dots,w_{16}$.}

Analogously we construct the cycles $w_9,\dots,w_{16}$ that correspond to submanifolds lying in the hyperplane
$x_4 = \frac{1}{2}$.
In this case the reasonings go through verbatim and it only needs to increase the numerical indices for
$w$-, $[L]$- and $[S]$-cycles by $8$ and to preserve them for $[T]$-cycles.

{\sc 3) The cycles $w_{17},\dots,w_{22}$.}

Let us define the following cycles:
$$
w_{17} = [T_{12}], \ \ \ w_{19} = [T_{13}], \ \ \ w_{21} = [T_{23}],
$$
$$
w^\prime_{18} = [S_{129(10)}] + [L_2] + [L_{10}] = -\frac{1}{2}\left( [L_1]-[L_2]+[L_9]-[L_{10}]\right) + \frac{1}{2} [T_{34}],
$$
$$
w^\prime_{20} = [S_{719(15)}] = - \frac{1}{2}\left([L_1]+[L_7]+[L_9]+[L_{15}]\right) - \frac{1}{2}[T_{24}],
$$
$$
w^\prime_{22}=[S_{13(11)9}]  = - \frac{1}{2}\left([L_1]+[L_3]+[L_9]+[L_{11}]\right) + \frac{1}{2} [T_{14}].
$$

We denote by $\Lambda_1$ the lattice generated by $w_1,\dots,w_{16}$, and by $\Lambda_2$ the lattice generated by
$w_{17},w^\prime_{18},w_{19},w^\prime_{20},w_{21},w^\prime_{22}$.

It follows from (\ref{form1}) and (\ref{form2}) that

1) $\Lambda_1$ and $\Lambda_2$ are pairwise orthogonal with respect to the intersection form:
$$
u \cap v = 0 \ \ \ \mbox{for all  $u \in \Lambda_1, v \in \Lambda_2$};
$$

2) in the basis $w_1,\dots,w_{16}$ the restriction of the intersection form onto $\Lambda_1$ looks like 
$$
\left(\begin{array} {cc} E_8(-1) & 0 \\ 0 & E_8(-1) \end{array}\right),
$$
where $E_8(-1)$ is given by (\ref{e8});

3) in the basis  $w_{17},w^\prime_{18},w_{19},w^\prime_{20},w_{21},w^\prime_{22}$ the restriction of the intersection form onto
$\Lambda_2$ is as follows
$$
\left(\begin{array}{cccccc} 0 & 1 & 0 & 0 & 0 & 0 \\  1 & -2 & 0 & -1 & 0 & -1 \\ 0 & 0 & 0 & 1 & 0 & 0 \\
0 & -1 & 1 & -2 & 0 & -1 \\ 0 & 0 & 0 & 0 & 0 & 1 \\ 0 & -1 & 0 & -1 & 1 & -2 \end{array}\right).
$$
It takes the form
$H \oplus H \oplus H$  in the basis $w_{19},\dots,w_{22}$, where the cycles $w_{17}$, $w_{19}$, and $w_{21}$ are the same
as before and 
$$
w_{18} = w^\prime_{18} + w_{17}, \ \ w_{20} = w^\prime_{20} + w_{17}+w_{19},  \ \ 
w_{22} = w^\prime_{22} + w_{17} + w_{19} + w_{21}.
$$

By definition,  the cycles $w_1,\dots,w_{22}$ are realized as formal sums of smooth sub manifolds and belong to 
$H_2(X;\Z)$. Since the intersection form on the lattice $\Lambda$ generated by these cycles takes the form (\ref{form}),
it is, in particular, is unimodular and therefore $\Lambda$ coincides with $H_2(X;\Z)$.

Therewith we proved

\begin{theorem}
In the basis $w_1,\dots,w_{22}$ the in\-tersection form
on $H_2(X;\Z)$
takes the canonical form (\ref{form}).
\end{theorem}

\vskip5mm

{\sc Remarks.} 1) In the definition of $S_{ijkm}$ we transform the tori $T_{ijkm}$ into
$T^\prime_{ijkm}$ by surgeries near fixed points of $\iota$  to achieve a situation when the intersection indices
of $S_{ijkm}$ with $L_n, n \in \{i,j,k,n\}$, are equal to one.
The homotopy $\varkappa_t$, which is used in this surgery, may be replaced by a similar one for which the resulted contour
$\varkappa_1$ bounds a negatively-oriented disc $D_0(L)$ in the plane $z_2 = \mathrm{const}$.
In this case the intersection index of the desingularized submanifold $S$ with $L$ would be equal to
$[S] \cap [L] = -1$ and
$$
[S] = [S_{ijkm}] + [L]
$$
in the homology group.
But in this case $[S]$ is not realized by complex submanifolds because it has 
a negative intersection index with the complex projective line $L$.
Moreover we can take another orientation of the sub manifold $S_{ijkm}$. This implies that if 
$$
[S_{ijkm}] = - \frac{1}{2}\left([L_i]+[L_j]+[L_k]+[L_m]\right) + \frac{1}{2}[T],
$$
where  $[T]$  is an integer combiaa=nation of cycles of the form $[T_{pq}]$,
then every cycle of the form 
$$
[S] = - \frac{1}{2}\left(\varepsilon_i [L_i] + \varepsilon_j [L_j] + \varepsilon_k [L_k] +
\varepsilon_m [L_m]\right) + \frac{1}{2}\varepsilon_ t [T],
$$
with $\varepsilon_i,\varepsilon_j,\varepsilon_k,\varepsilon_m,\varepsilon_t \in \{1, -1\}$,
is realized by a sub manifold which is homeomorphic to the sphere.

2) By the Hurewicz theorem, for simply connected manifolds the natural homomorphism
$$
\pi_2(X) \to H_2(X;\Z)
$$
is an isomorphism. By the previous remark, the cycles $w_1,\dots,w_{16}$ belong to the image of this homomorphism. The formula (\ref{s}) demonstrates how the cycles of the form $[T_{pq}]$ are represented by sums of spherical cycles. For instance, it follows from
$$
[S_{1357}] = - \frac{1}{2}\left([L_1]+[L_3]+[L_5]+[L_7]\right)  +\frac{1}{2}[T_{12}]
$$
that
$$
[T_{12}]=  2[S_{1357}] + [L_1]+[L_3]+[L_5]+[L_7].
$$

3) For a closed oriented $N$-dimensional smooth manifold $X$ the intersection of cycles
$$
H_k (X;\Z) \times H_l(X;\Z) \stackrel{\cap}{\longrightarrow} H_{k+l-N}(X;\Z)
$$
is dual to the cohomological product in the following sense. Indeed, if the cycles
$u \in H_k(X;\Z)$ and $v \in H_l(X;\Z)$ are realized by smooth submanifolds $Y$ and $Z$ that intersect each other transversally, then their intersection is a smooth submanifold $W$ of dimension  $k+l-N$ and
$$
Du \cup Dv = Dw,
$$
where $w$ is the cycle realized by $W$ and
$$
D: H_i(X;\Z) \to H^{N-i}(X;\Z), \ \ \ i=0,\dots,N,
$$
is the Poincare duality. The same holds for nonoriented manifolds for the homo\-lo\-gy with
$\Z_2$ coefficients. It is known that not all cycles are realized by smooth submanifolds; however this duality
is generalized for all manifolds by using special types of chains. An accurate and complete exposition of
this const\-ruc\-tion is given \cite{GP}, where the homology group $H_\ast(X)$ with the intersection operation is called the (intersection) Lefschetz ring and, in par\-ti\-cu\-lar, by using an isomorphism $D$ of it to $H^\ast(X)$ the topological invariance
of the intersection ring is established.

4) For an eight-dimensional Kummer manifold, which is ob\-tained by a desingularization of
the quotient-space $T^8/\Z_2$, the analogs of the cycles $S_{ijkm}$ of the half dimension of the manifold are four-cycles
generated by tori $T^4 \subset T^8$. These four-dimensional tori pass through 16 fixed points of $\iota$ and are desingularized into submanifolds diffeomorphic to $K3$ surfaces. Analogously to the case of $K3$ surfaces (see Remarks 1 and 2) we may conclude that all four-dimensional cycles are realized by linear combinations of embedded $K3$ surfaces and complex projective planes $\C P^2$.

\end{document}